\documentclass[12pt]{slides}
\usepackage{amsmath,amsfonts,amssymb,amsthm}
\begin{document}
\author{{\bf Denis~I.~Saveliev\/}}

\title{{\bf Choice~and~Regularity:\\
Common~Consequences in~Logic\/}}
\date{21 August 2007, Xi'an
\vfill{\begin{tiny}
Partially supported by grant~06-01-00608-a 
of Russian Foundation for Basic Research
\end{tiny}
}}
\maketitle

\theoremstyle{plain}
\newtheorem{thm}{Theorem}
\newtheorem{prb}{Problem}
\newtheorem{coro}{Corollary}
\newtheorem*{tm}{Theorem}
\newtheorem*{cor}{Corollary}
\newtheorem*{lm}{Lemma}
\newtheorem*{fct}{Fact}
\newtheorem*{fcts}{Facts}

\theoremstyle{definition}
\newtheorem*{df}{Definition}
\newtheorem*{rfr}{References}

\theoremstyle{remark}
\newtheorem*{rmk}{Remark}
\newtheorem*{exm}{Example}

\newcommand{\eqs}{ {\;=^{\!*}\,} }
\newcommand{\lesss}{ {\;<^{\!*}\,} }
\newcommand{\leqs}{ {\;\le^{\!*}\,} }
\newcommand{\wo}{\mathop {\mathrm {wo\,}}\nolimits }
\newcommand{\pwo}{\mathop {\mathrm {pwo\,}}\nolimits }
\newcommand{\wf}{\mathop {\mathrm {wf\,}}\nolimits }
\newcommand{\ewf}{\mathop {\mathrm {ewf\,}}\nolimits }
\newcommand{\cf}{ {\mathop{\mathrm {cf\,}}\nolimits} }
\newcommand{\cl}{ {\mathop{\mathrm {cl\,}}\nolimits} }
\newcommand{\tc}{ {\mathop{\mathrm {tc\,}}\nolimits} }
\newcommand{\dom}{ {\mathop{\mathrm {dom\,}}\nolimits} }
\newcommand{\ran}{ {\mathop{\mathrm {ran\,}}\nolimits} }
\newcommand{\fld}{ {\mathop{\mathrm {fld\,}}\nolimits} }
\newcommand{\otp}{ {\mathop{\mathrm {otp\,}}\nolimits} }
\newcommand{\lh}{ {\mathop{\mathrm {lh\,}}\nolimits} }
\newcommand{\gt}{\mathfrak}
\newcommand{\ext}{\mathrm{ext}}
\newcommand{\rk}{ {\mathrm{rk}} }
\newcommand{\Z}{ {\mathrm Z} }
\newcommand{\ZF}{ {\mathrm {ZF}} }
\newcommand{\ZFA}{ {\mathrm {ZFA}} }
\newcommand{\ZFC}{ {\mathrm {ZFC}} }
\newcommand{\E}{ {\mathrm {AE}} }
\newcommand{\AR}{ {\mathrm {AR}} }
\newcommand{\GR}{ {\mathrm {GR}} }
\newcommand{\WR}{ {\mathrm {WR}} }
\newcommand{\WO}{ {\mathrm {WO}} }
\newcommand{\GWO}{ {\mathrm{GWO}} }
\newcommand{\LO}{ {\mathrm {LO}} }
\newcommand{\AF}{ {\mathrm {AF}} }
\newcommand{\AC}{ {\mathrm {AC}} }
\newcommand{\DC}{ {\mathrm {DC}} }
\newcommand{\GC}{ {\mathrm {GC}} }
\newcommand{\ACM}{ {\mathrm {ACM}} }
\newcommand{\GCM}{ {\mathrm {GCM}} }
\newcommand{\AD}{ {\mathrm {AD}} }
\newcommand{\AInf}{ {\mathrm {AInf}} }
\newcommand{\AU}{ {\mathrm {AU}} }
\newcommand{\AP}{ {\mathrm {AP}} }
\newcommand{\CH}{ {\mathrm {CH}} }
\newcommand{\APr}{ {\mathrm {APr}} }
\newcommand{\ASp}{ {\mathrm {ASp}} }
\newcommand{\ARp}{ {\mathrm {ARp}} }
\newcommand{\AFA}{ {\mathrm {AFA}} }
\newcommand{\BAFA}{ {\mathrm {BAFA}} }
\newcommand{\FAFA}{ {\mathrm {FAFA}} }
\newcommand{\SAFA}{ {\mathrm {SAFA}} }
\newcommand{\NS}{ {\mathrm{NS}} }
\newcommand{\Cov}{ {\mathrm{Cov}} }
\newcommand{\Con}{ {\mathrm{Con}} }
\newcommand{\PA}{ {\mathrm{PA}} }
\newcommand{\BF}{ {\mathrm{BF}} }
\newcommand{\RP}{ {\mathrm{RP}} }
\newcommand{\RC}{ {\mathrm{RC}} }
\newcommand{\RE}{ {\mathrm{RE}} }
\newcommand{\J}{ {\mathrm {J}} }
\newcommand{\Ja}{ {\mathrm {Ja}} }
\newcommand{\Jb}{ {\mathrm {Jb}} }
\newcommand{\Jc}{ {\mathrm {Jc}} }
\newcommand{\Sc}{ {\mathrm {Sc}} }
\newcommand{\wSc}{ {\mathrm {wSc}} }
\newcommand{\Sk}{ {\mathrm {Sk}} }
\newcommand{\GSc}{ {\mathrm {GSc}} }
\newcommand{\GwSc}{ {\mathrm {GwSc}} }
\newcommand{\GSk}{ {\mathrm {GSk}} }
\newcommand{\T}{ {\mathrm {T}} }
\newcommand{\TST}{ {\mathrm {TST}} }
\newcommand{\PI}{ {\mathrm {PI}} }
\newcommand{\LF}{ {\mathrm {LF}} }
\newcommand{\LR}{ {\mathrm {LR}} }

\newpage

It is well-known that
Choice and Regularity
are {\it independent\/} of each other but 
have important {\it common consequences\/}
of logical character
(reflection principles, 
representations of classes by sets, etc.).
In my talk,
I~shall try:

(A)
{\it 
To explain this phenomenon\/},

(B)
{\it
To consider relationships between these consequences 
(and near principles) in detail\/},

and besides,

(C)
{\it
To consider some arguments related to 
truth of various principles in set theory\/}.

All theorems can be proved 
in~$\ZF$,
the Zermelo-Fr\"ankel set theory,
{\it minus\/} Regularity
(quite often in some its fragments).

\newpage
{\centerline{\Large{\bf Choice and Regularity\/}}}

\newpage
Basic definitions:

A~function
$F:X\to\bigcup X$ is
a~{\it choice function\/} iff
$F(x)\in x$
for all nonempty $x\in X$.

A~relation $R\subseteq X\times X$ is
{\it well-founded\/} iff
each nonempty subset $S\subseteq X$ has
an $R$-minimal element
(i.e., $x\in S$ such that $\neg(y\,R\:x)$ 
for all $y\in S-\{x\}$).

An ordering~$\le$ is 
a~{\it well-ordering\/} of~$X$ iff
each nonempty subset $S\subseteq X$ has
a~$\le$-least element
(i.e., $x\in S$ such that $x\le y$ 
for all $y\in S$).

Equivalently, 
$\le$ is linear and well-founded.

\newpage

\noindent
{\bf The Axiom of Choice, AC.\/}
{\it For any set
there is a~choice function on~it.\/}

There are a~number of equivalent principles,
the most famous are perhaps
Zorn's Lemma (Kuratowski) and

{\bf The Well-Ordering Principle, WO\/}.
{\it For any set there is a~well-ordering on~it.\/}

\noindent
\begin{tm}[Zermelo]
$\AC$ is equivalent to~$\WO$.
\end{tm}

AC has a~deep impact on the universe of set theory
by giving as nice consequences, e.g.,
$$
\text{{\it All cardinals form a~well-ordered hierarchy\/}}
$$
as well as ugly ones, e.g.,
$$
\text{The Banach-Tarski Paradox.}
$$

\newpage

{\bf The Axiom of Regularity, AR.\/}
{\it Any nonempty set has 
an $\in$-minimal element.\/}

A~set is {\it well-founded\/} iff
$\in$ is well-founded on its transitive closure.
Then AR is equivalent to the sentence:
$$
\text{{\it All sets are well-founded}\/}
$$
(and another name of~AR is 
the {\it Axiom of Foundation\/}).

To formulate an equivalent principle,
recall the {\it cumulative hierarchy\/} of sets:
\begin{align*}
V_0
&=
\emptyset,
\\
V_{\alpha+1}
&=
P(V_\alpha),
\\
V_\alpha
&=
\bigcup_{\beta<\alpha}V_\beta
\text{ if $\alpha$ is limit.}
\end{align*}

\noindent
\begin{tm}[von Neumann]
$\AR$ is equivalent to
$$
V=\bigcup_{\alpha\in Ord}V_\alpha.
$$
\end{tm}

\newpage

If 
$R\subseteq X\times X$ 
is well-founded,
we have:
\\
\noindent
(i)
$R$-Induction
\\
\noindent
(ii)
$R$-Recursion
\\
\noindent
(iii)
The~{\it rank\/} function
$$
\rk_R:(X,R)\to(Ord,<)
$$
i.e., 
a~strong homomorphism
stratifying~$X$ into levels~$X_\alpha$:
$$
X=\bigcup_{\alpha}X_\alpha
$$
\\
\noindent
(iv)
The {\it transitive collapse\/}
$$
\pi_R:
(X,R)\to(\bigcup_{\alpha}V_\alpha,\in)
$$
which is a~strong homomorphism
allowing to get

\noindent
\begin{tm}[The Mostowski Collapsing Lemma]
Any extensional well-founded relation
is isomorphic to a~unique transitive one
(and a~unique possible isomorphism is its transitive collapse).
\end{tm}

So, AR has mainly 
a~``simplifying'' character:
we can use all these nice properties.

\newpage

Quite often principles have
local and global forms.
Typically,
a~{\it local\/}/{\it global\/} principle
says about
sets/classes
or
one formula/all formulas.
Global versions of the previous principles:

{\bf The Global Choice, GC.\/}
{\it There is a~choice function on the universe.\/}

{\bf The Global Well-Ordering, GWO.\/}
{\it There is a~well-ordering of the universe.\/}

{\bf The Global Regularity, GR.\/}
{\it Any nonempty class has 
an $\in$-minimal element.\/}

(GR~is a~schema.
To formulate~GC,
we add a~new functional symbol.
For~GWO,
we add a~new predicate symbol.)

\newpage

{\bf Lemma.\/}
{\it
\\
1.
$\GR$ is equivalent to~$\AR$.
\\
2.
$\GWO$ implies $\GC$.
\\
3.
$\GC$ implies $\AC$.
\\
4.
$\AC+\neg\GC$ is consistent.
\\
5.
$\GC+\AR$ implies $\GWO$.
\\
6.
$\GWO+\neg\AR$ is consistent.
\\
7.
$\neg\AC+\AR$ is consistent.
\\
8.
$\neg\AC+\neg\AR$ is consistent.
\/}

(The only hard clause is~(4).
Later I~shall show that
one can sharp~(5)
by replacing ``implies'' with ``is equivalent to''
and AR~with a~weaker principle~BF.)

\newpage

To complete this account,
note that
Choice {\it plus\/} Regularity together
can be formulated in a~single way:

{\bf The Choice of Minimals, ACM.\/}
{\it For any set~$X$
there is a~choice function~$F$ on~$X$
such that
$F(x)\cap x=\emptyset$
for all nonempty $x\in X$.\/}

{\bf The Global Choice of Minimals, GCM.\/}
{\it \\There is a~choice function~$C$ on~$V$
such that
\\\text{$C(x)\cap x=\emptyset$}
for all nonempty sets~$x$.\/}

Clearly,

{\bf Lemma.\/}
{\it
\\
1.
$\ACM$ is equivalent to $\AC+\AR$.
\\
2.
$\GCM$ is equivalent to $\GC+\AR$. 
\/}

\newpage
{\centerline{\Large{\bf Best-Foundedness\/}}}

\newpage

To explicate 
why Choice (mainly in the stron\-gest form~GWO)
and Regularity have common consequences,
I~isolate their ``intersection'':
a principle (called here Best-Foundedness)
which is consistent with
negations of both axioms 
but implies all these consequences.

\newpage

Let me say that
a~well-founded relation~$E$ is
{\it best-founded\/} iff
$\{x:\rk_E(x)=\alpha\}$
is a~set
for every ordinal~$\alpha$.

By Replacement, 
then  
$U_\alpha=\{x:\rk_E(x)<\alpha\}$ 
is also a~set
for every~$\alpha$.

{\it Examples.\/}
The empty relation on a~proper class
is well- but not best-founded.
$\in$ is best-founded on transitive well-founded sets,
and so (by the Mostowski theorem)
all extensional well-founded relations are best-foun\-ded.

{\bf The Best-Foundedness Axiom, BF.\/}
{\it There is a~best-founded relation on~$V$.\/}

(The axiom is in the language 
with a~new predicate symbol.)

{\bf Lemma.\/}
{\it
\\
1.
$\AR$ implies $\BF$.
\\
2.
$\GWO$ implies $\BF$.
\\
3.
$\BF+\neg\AC+\neg\AR$ is consistent.
\/}

\newpage

The principle has a~number of reformulations
(in appropriate languages).
Define:

$A$ is
{\it club\/} iff
it is
$\subseteq$-cofinal in~$V$ and
for any $\subseteq$-directed $x\subseteq A$
we have $\bigcup x\in A$.

$A$ is a~{\it basis\/}
iff
$\{P^\alpha(x):x\in A\,\wedge\,\alpha\in Ord\}$
is $\in$-cofinal in~$V$.

{\it Example.\/}
$\AR$ is equivalent to any of (1) and~(2):
\\
1.
$\{V_\alpha:\alpha\in Ord\}$
is club.
\\
2.
$V_1$ ($=\{\emptyset\}$) is a~basis.

Under~BF
the sets $U_\alpha$ play 
much the same part 
as the sets $V_\alpha$ under~AR.
E.g.,
$\{U_\alpha:\alpha\in Ord\}$
is club.

\newpage

Moreover,

{\bf Lemma.}
{\it
$\BF$ is equivalent to any of (1)--(4):
\\
1.
There is a~function
$F:V\to Ord$
such that
$F^{-1}(\alpha )$ is a~set
for all~$\alpha$.
\\
2.
There is a~well-ordered 
$\in$-cofinal in~$V$ class.
\\
3.
There is a~well-ordered club class.
\\
4.
There is a~well-ordered basis.
\\
5.
There is a~well-ordered partition 
of~$V$ into sets.
\/}

Clause~(2) gives
a~visual notion about~BF:
intuitively,
ordinals of a~model show its ``height'';
then 
a~model witnessing~BF is ``stretched upward'' 
while
a~model refuting~BF is ``inflated in width''.

\newpage

A~similarity:
GWO well-orders the whole universe
while BF~well-orders some its ``essential'' part
(a~basis or a~club).
Moreover,
this can be maked
in a~natural way:

{\bf Lemma.\/}
{\it
If
there is a~well-orderable class
that is club (or a~basis),\,
then
there is such a~class
which is moreover
$\in$-~and $\subseteq$-well-ordered.
\/}

There are less obvious reformulations of~BF,
one of which (concerning the ordinal definability)
I~shall give a~bit later.

Finally,
BF is exactly
what is missing in~GC to be~GWO:

{\bf Theorem.}
{\it
$\GC+\BF$ is equivalent to $\GWO$.\/}

(Cf.~with (5) of Lemma above.)

\newpage

{\centerline{\Large{\bf Consequences\/}}}

\newpage

Showing that BF~works,
I~shall consider 
following its consequences:

The existence of Skolem and Scott functions,
\\
The reflection of formulas at sets,
\\
The expressibility of the ordinal definability, 
\\
The representability of equivalence classes by sets,

and relationships between them.

\newpage

Let
$\varphi (u,\ldots,x)$
be a~formula with the parameters
$u,\ldots,x$.

A~function $A_\varphi$
is a~{\it Skolem function for~$\varphi$\/} 
iff
$$
(\exists x)\;\varphi(u,\ldots,x)
\;\to\;
\varphi(u,\ldots,A_\varphi(u,\ldots)).
$$

Similarly,
let me say that
a~function $B_\varphi$
is a {\it Scott function for~$\varphi$\/} 
iff
$$
(\exists x)\;\varphi(u,\ldots,x)
\;\to\;
(\exists x\in B_\varphi(u,\ldots))\;
\varphi(u,\ldots,x)
$$
and
$$
(\forall x\in B_\varphi (u,\ldots))\;
\varphi(u,\ldots,x).
$$

\newpage

Thus
$A_\varphi$ chooses
a~single point from the class
$\{x:\varphi(u,\ldots,x)\}$:
$$
A_\varphi(u,\ldots)
\,\in\,
\{x:\varphi(u,\ldots,x)\}
$$
while
$B_\varphi$ separates from this class
its subset
$$
B_\varphi(u,\ldots)
\,\subseteq\,
\{x:\varphi(u,\ldots,x)\}
$$
such that
the set
$B_\varphi(u,\ldots)$ is nonempty
whenever
the class
$\{x:\varphi(u,\ldots,x)\}$ so~is.

{\it Remark.\/}
Scott was probably first who noted
that such functions can be used instead of
Skolem functions in absence of~AC.

\newpage

Consider the following schemas
(in extended languages):

{\bf The Skolem Principle, Sk.\/}
{\it For any formula
there is a~Skolem function.\/}

{\bf The Scott Principle, Sc.\/}
{\it For any formula
there is a~Scott function.\/}

$\Sk_\varphi$ and $\Sc_\varphi$ denote 
the instances of these schemas.

{\bf Lemma.}
{\it
\\
1.
$\Sk_\varphi$ implies~$\Sc_\varphi$.
\\
2.
$\AC+\Sc_\varphi$ implies~$\Sk_\varphi$.
\\
3.
$\BF$ implies~$\Sc$.
\\
4.
$\AC+\BF$ implies~$\Sk$.
\/}

\newpage

Via coding,
one can formulate global variants
of Skolem and Scott functions
(uniformly for all formulas).
A~global Skolem function 
acts like 
a~choice function on definable classes
while
a~global Scott function  
separates subsets from them.
Let GSk and~GSc denote 
the global variants of Sk and~Sc.

{\bf Lemma.}
{\it
\\
1.
$\GSk$ is equivalent to~$\GC+\GSc$.
\\
2.
$\GSk+\BF$ is equivalent to~$\GWO$.
\/}

\newpage

I~need Scott (or Skolem) functions mainly
to have Reflection.

A~class~$M$ {\it reflects\/}
a~formula $\varphi(x,\ldots)$
iff
for all $x,\ldots\in M$
$$
\varphi^M(x,\ldots)
\;\leftrightarrow\;
\varphi(x,\ldots).
$$

{\bf The Reflection Principle, RP.\/}
{\it Each formula is reflected at some set.\/}

(RP is a~schema, 
$\RP_\varphi$ are instances.)

It follows from RP that
each true formula
has a~set model.

\newpage

Of course,
the principle holds for 
finitely many formulas as well.
Let me rewrite it as follows:
If $\Gamma$ is a~finite set of formulas,
then there is a~set~$M$ such that
$$
M\;\prec_\Gamma\;V.
$$

Thus
RP is a~local variant of
the L\"owenheim-Skolem Theorem.
But unlike it,
RP~can be proved {\it inside\/} (some) set theory:

{\bf Theorem.\/}
{\it
$\Sc$ implies $\RP$.\/}

(Take a~{\it Scott hull\/}.)

Moreover, 
$\Sc$ gives a~club class 
of reflecting sets, and
$\BF$ gives a~club class 
of reflecting sets of form~$U_\alpha$.

{\it Remarks.\/}
1.
Without Choice,
we know nothing about 
the size of submodels.
\\
2.
The full L\"owenheim-Skolem Theorem
(without an evaluation of the cardinality)
can be obtained in the same way
as a~metatheorem.

\newpage

We consider two applications of Reflection.
The first concerns the finite axiomatizability:

Let us call a~theory 
{\it sufficiently rich\/} iff
it admits a~coding.
(E.g., 
ZF minus Infinity so is).

{\bf Proposition.}
{\it
Let $\T$ be
sufficiently rich
consistent theory and
$\T\vdash\RP$.
Then
$\T$ is not finitely axiomatizable.
\/}

(Apply the Second Incompleteness Theorem.)

{\it Examples.\/}
The theory consisting of
Union, Power Set, Replacement, 
and Best-Founded\-ness
is not finitely axiomatizable.
The same for any its consistent extension
(e.g.,~ZF).
On the other hand,
in the Zermelo set theory~Z 
(which is finitely axiomatizable),
RP is not provable.

\newpage

Another application of~RP:
the description of ordinal-definable sets
inside of set theory.

A~class is
{\it ordinal-definable\/} iff
it is of the form
$\{u:\varphi(u,\alpha,\ldots)\}$
for some formula~$\varphi$
where all $\alpha,\ldots$ are ordinals.
$OD$ is the class of all ordinal-definable sets.
$\cl$ is the closure
under G\"odel operations.

A~well-known fact:
$\AR$ implies
$$
OD=\cl(\{V_\alpha:\alpha\in Ord\}).
$$
It follows that
$OD$~is well-orderable and club
(and moreover,
the largest inner model of~ZF
with a~global well-ordering
definable via~$\in$).

\newpage

We sharp:

{\bf Theorem.\/}
{\it
$\BF$~implies\/}
$$
OD=\cl(\{U_\alpha:\alpha\in Ord\}).
$$

(Use~RP to prove $\subseteq$.)

{\bf Corollary.\/}
{\it
$\BF$ holds
\,iff\,
$OD$~is well-orderable and club.
\/}

Thus again 
(like the characteristic of GWO via BF and~GC)
we sharp an old result of form
$$
\Gamma+\AR\,\text{ {\it implies\/} }\,\Delta
$$
by a~new result of form
$$
\Gamma+\BF\,\text{ {\it is equivalent to\/} }\,\Delta
$$
(where $\Gamma$ and $\Delta$ are some sets of formulas).
This supports a~naturality of~BF.

\newpage

As the last interesting consequence of~$\BF$,
consider representations of equivalence classes by sets.

Let $\varphi (x,y)$
define an equivalence:
$$
\varphi(x,y)
\wedge
\varphi(y,z)
\;\to\;
\varphi(y,x)
\wedge
\varphi(x,z).
$$
A~function~$F_\varphi$
{\it represents\/} the equivalence 
defined by~$\varphi$ iff
$$
\varphi(x,y)
\:\leftrightarrow\;
F_\varphi(x)=F_\varphi(y).
$$

{\bf The Representation of Classes Principle, RC.\/}
{\it For any equivalence formula
there is a~representing function.\/}

(RC is a~schema, 
$\RC_\varphi$ are instances.)

\newpage

Of course,
$$
\Sc_\varphi\,\text{ implies }\,\RC_\varphi
$$
since any Scott function for~$\varphi$ 
represents the equivalence
in a~``natural way''.
But unlike Scott functions,
$F_\varphi(x)$ does not meet necessarily
the equivalence class 
$\{y:\varphi(x,y)\}$.

Sometimes
$\AC$ suffices for some instances of~RC:

{\it Example.\/}
If $\varphi$ expresses the same cardinality,
then $\AC$ implies $\Sk_\varphi$ 
and so~$\RC_\varphi$.

Moreover,

{\bf Theorem.\/}
{\it
Let $\varphi$ define an equivalence.
Then $\AC+\RP_\varphi$ implies~$\Sk_\varphi$.
\/}

{\bf Corollary.\/}
{\it
$\AC+\RP$ implies~$\RC$.\/}

\newpage

{\bf Question.\/}
{\it
Is any of the following implications provable:
\\
1.
$\AC$ implies $\Sc$?
\\
2.
$\AC+\Sc$ implies $\Sk$?
\\
3.
$\Sc$ implies $\BF$?
\\
4.
$\GC$ implies $\BF$? 
\\
5.
$\GC$ implies $\RC$?
if No for (1)--(4).
\/}

{\bf Conjecture.\/}
{\it 
No for all (1)--(5).\/}

(Partial results.)

To prove consistency results of such kind,
I~develop a~method of construction of models
via {\it automorphism filters\/}
(a~generalization of well-known permutation model method).
Two main technical obstacles
arising without Regularity
(or with a~proper class of atoms):
\\
(i)
Replacement,
\\
(ii)
The Transfer Theorem (Jech--Sochor).

\newpage
{\centerline{\Large{\bf Truth and\/}}} 

\noindent
{\centerline{\Large{\bf Interpretability Strength\/}}}

\newpage

ZFC is highly incomplete:
there are very natural set-theoretical questions
independent of it,
the most famous of which is perhaps
the Continuum Hypothesis
(which is really
a~{\it third order arithmetical\/} sentence).

To complete ZFC,
it was proposed 
a~number of principles of various kind
having important consequences:
e.g.,

Large Cardinal Axioms,
\\
The Axiom of Determinacy (Mycielski),
\\
Proper Forcing Axioms (Shelah),
\\
The $\Omega$ Conjecture (Woodin),
\\
Generic Large Cardinals (Foreman),
\\
The Inner Model Hypothesis (Sy Friedman),
\\
etc.

Is there a~general criterion
for rejecting/accepting 
such a~principle
as a~{\it true\/} axiom 
about {\it all\/} sets?

\newpage

Let me point out
a~simple criterion 
indicating some principles as
{\it surely wrong\/}.
An idea:
since all mathematical objects are sets,
an ideal TST (``True Set Theory'') 
must capture all possible theories.
Hence having some $\T\subseteq\TST$
and examining a~new principle $\Gamma$,
we must reject $\Gamma$
if it restricts this possibility:

{\it If\, 
$\T+\,\Gamma$ loses 
the interpretability strength of $\T$
then $\Gamma$ is wrong.\/}

Here:
An extension $\T_1\supseteq\T$ 
of a~theory~$\T$
{\it loses the interpretability strength\/} of~$\T$
iff
there is $\T_2\supseteq\T$
non-interpretable in
all $\T_3\supseteq\T_1$.
Otherwise
$\T_1$ extends~$\T$ 
{\it without loss of the interpretability strength\/}.

\newpage

{\it Example.\/}
``All sets are constructible''
$$V=L$$
is a~nice axiom
since it looks ``empirically complete'',
but it loses the interpretability strength
of~ZFC:
Under $V=L$,
there is no measurable cardinals,
even in inner models.
So, it is wrong.
May be so is 
``All sets are ordinal-definable''
$$V=\,OD\,?$$
or even
``All sets are well-founded''
$$V=\,\bigcup_\alpha V_\alpha\,?$$
(See Question below.)

\newpage

To find a~criterion 
indicating some principles as
{\it surely true\/}
is much more hard.
(Of course,
if $\varphi$ is surely wrong
then $\neg\varphi$ is surely true
but too {\it noneffective\/} as a~rule.
E.g.,
cf.~``there exists a~nonconstructible set''
with
``\,$0^\sharp$~exists''.)
Advancing in the same way,
we can describe only 
{\it possibly true\/} principles:

{\it If\, 
$\T+\Gamma$ extends~$\T$ 
without loss of the interpretability strength
and does not interpretable in~$\T$, 
then $\Gamma$ can be true\/}.

Such candidates for being true
cannot be really true all together
because contradict to each other.
But some of them are concordant:

{\it Example.\/}
Current large cardinal axioms form
a~well-ordered hierarchy.
(An empirical fact;
Woodin offers a~partial explication.)

\newpage

{\it Remark.\/}
Criterions based on
interpretability are perhaps 
most important but not sufficient.
E.g.,
put
$\T$ be
$$
\ZF-\text{Infinity}+\neg\,\text{Infinity}+\Con(\ZFC).
$$
$\ZFC$ is interpretable in~T,
but I~think
this theory
(in which infinite objects do not exist)
is not a~correct theory of {\it all\/} sets.
Likewise for extensions of~$\ZFC$
by large cardinals.
(An explication is beyond my talk.)

\newpage

Among the axioms of~ZFC,
only Extensionality~(AE) and Regularity
have an ``impoverishing'' character
(since forbid sets of certain structure);
a~character of~AC is unclear;
and other axioms have 
an ``enriching'' character
(since permit to construct new sets).
Does this ``impoverishment''
really decrease
the interpretability strength?
Under GWO,
we have easily the answer No:

{\bf Lemma.\/}
{\it
$\ZF-\E-\AR+\GWO$ can be exten\-ded 
by $\E+\AR$
without loss of 
the interpretability strength.\/}

{\bf Question.\/}
{\it
Can one extend 
without loss of
the interpretability strength:
\\
1.
$\ZF-\E-\AR$
by $\BF$?
\\
2.
$\ZF-\E-\AR+\BF$ by $\E+\AR$?
\\
3.
$\ZF$ by $\GWO$?
\/}

\end{document}